\documentclass{article}
\usepackage{stmaryrd}
\usepackage{mathrsfs}
\usepackage[centertags]{amsmath}
\usepackage{amsfonts, dsfont}
\usepackage{amssymb}
\usepackage{amsthm}

\input xy
\xyoption{all}

\theoremstyle{plain}
\newtheorem{thm}{Theorem}[section]
\newtheorem{cor}[thm]{Corollary}
\newtheorem{lem}[thm]{Lemma}
\newtheorem{prop}[thm]{Proposition}
\newtheorem*{thm2}{Theorem}
\newtheorem*{cor2}{Corollary}

\theoremstyle{definition}
\newtheorem{defn}[thm]{Definition}
\newtheorem*{remark}{Remarks}
\newtheorem*{ack}{Acknowledgments}

\newcommand{\bd}{\begin{defn}}
\newcommand{\ed}{\end{defn}}
\newcommand{\bl}{\begin{lem}}
\newcommand{\el}{\end{lem}}
\newcommand{\bp}{\begin{prop}}
\newcommand{\ep}{\end{prop}}
\newcommand{\bt}{\begin{thm}}
\newcommand{\et}{\end{thm}}
\newcommand{\bc}{\begin{cor}}
\newcommand{\ec}{\end{cor}}
\newcommand{\br}{\begin{remark}}
\newcommand{\er}{\end{remark}}
\newcommand{\bdi}{\begin{diagram}}
\newcommand{\edi}{\end{diagram}}
\newcommand{\beq}{\begin{equation}}
\newcommand{\eeq}{\end{equation}}
\newcommand{\ba}{\begin{array}}
\newcommand{\ea}{\end{array}}
\newcommand{\bpf}{\begin{proof}}
\newcommand{\epf}{\end{proof}}

\newcommand{\Z}{\mathds{Z}}
\newcommand{\Q}{\mathds{Q}}
\newcommand{\Zp}{\mathds{Z}_{p}}
\newcommand{\Qp}{\mathds{Q}_{p}}

\newcommand{\Ep}{E[p^{\infty}]}
\newcommand{\al}{\alpha}
\newcommand{\be}{\beta}
\newcommand{\Ga}{\Gamma}

\newcommand{\La}{\Lambda}
\newcommand{\la}{\lambda}
\newcommand{\e}{\varepsilon}
\newcommand{\si}{\sigma}
\newcommand{\Si}{\Sigma}
\newcommand{\s}{\overrightarrow{s}}
\newcommand{\mF}{\mathcal{F}}

\DeclareMathOperator{\Sel}{Sel} \DeclareMathOperator{\Gal}{Gal}
\DeclareMathOperator{\Hom}{Hom} \DeclareMathOperator{\rank}{rank}
\DeclareMathOperator{\corank}{corank}
\DeclareMathOperator{\Char}{Char}

\newcommand{\cts}{\mathrm{cts}}

\newcommand{\ot}{\otimes}
\newcommand{\ilim}{\displaystyle \mathop{\varinjlim}\limits}

\newcommand{\im}{\mathrm{im}\,}

\newcommand{\lra}{\longrightarrow}

\newcommand{\ps}[1]{\llbracket #1 \rrbracket}

\begin{document}
\title{On the algebraic functional equation of the eigenspaces of mixed signed Selmer groups of elliptic curves with good reduction at primes above $p$}
 \author{ Suman Ahmed\footnote{School of Mathematics and Statistics,
Central China Normal University, Wuhan, 430079, P.R.China.
 E-mail: \texttt{npur.suman@gmail.com}}  \quad
  Meng Fai Lim\footnote{School of Mathematics and Statistics $\&$ Hubei Key Laboratory of Mathematical Sciences,
Central China Normal University, Wuhan, 430079, P.R.China.
 E-mail: \texttt{limmf@mail.ccnu.edu.cn}} }
\date{}
\maketitle

\begin{abstract} \footnotesize
\noindent Let $p$ be an odd prime number, and let $E$ be an elliptic curve defined over a number field which has good reduction at every prime above $p$.
Under suitable assumptions, we prove that the $\eta$-eigenspace and the $\bar{\eta}$-eigenspace of mixed signed Selmer group of the elliptic curve are pseudo-isomorphic. As a corollary, we show that the $\eta$-eigenspace is trivial if and only if the $\bar{\eta}$-eigenspace is trivial. Our results can be thought as a reflection principle which relate an Iwasawa module in a given eigenspace
with another Iwasawa module in a ``reflected'' eigenspace.


\end{abstract}

\section{Introduction}

The main conjecture of Iwasawa theory predicts a relation between a
Selmer group and a conjectural $p$-adic $L$-function (see \cite{BL, G89, K, Kob, SU}). This $p$-adic
$L$-function is expected to satisfy a conjectural functional equation in a certain sense.
In view of the main conjecture and this conjectural functional equation, one would
expect to have certain algebraic relationship between the corresponding Selmer groups. In this paper, we will examine this phenomenon for the eigenspaces of the mixed signed Selmer groups of elliptic curves with good reduction at primes above $p$. We shall describe our main result briefly in this introductory section.

Throughout the paper, $p$ will denote a fixed odd prime. Let $F'$ be a number field and $E$ an elliptic curve defined over $F'$. Let $F$ be a finite extension of $F'$. The following assumptions will be in full force throughout the paper.

\begin{itemize}
\item[(S1)] The elliptic curve $E$ has good reduction at all primes of $F'$ above $p$, and at least one of which is a supersingular reduction prime of $E$.

 \item[(S2)] For each prime $u$ of $F'$ above $p$ at which the elliptic curve $E$ has supersingular reduction, we assume that the following statements are valid.

 (a) $F'_u=\Qp$.

 (b) $a_u = 1 + p - |\tilde{E}_u(\mathbb{F}_p)| = 0$, where $\tilde{E}_u$ is the reduction of $E$ at $u$.

 (c) $u$ is unramified in $F/F'$.
\end{itemize}

Denote by $\mF_n = F(\mu_{p^{n+1}})$ and $\mF_{\infty} = \cup_n \mF_n$. We define a signed Selmer group $\Sel^{\s}(E/\mF_{\infty})$ (see the body of the paper for the precise definition) and write $X^{\s}(E/\mF_{\infty})$ for its Pontryagin dual. We shall write $\Si^{ss}_p$ for the primes of $F$ above $p$ at which $E$ has good supersingular reduction. For a character $\eta$  of $\Gal(F(\mu_p)/F)$, we write $e_{\eta}$ for the idempotent $\frac{1}{|G|}\sum_{\si\in G} \eta(\si)\si^{-1}$ and $\bar{\eta}$ for its contragradient. The goal of the paper is to prove the following result.

\begin{thm2}[Theorem \ref{main theorem}]
Let $\eta$ be a character of $\Gal(F(\mu_p)/F)$. For each $v\in \Si_p^{ss}$ with $s_v=+$, assume that $|F_v:\Qp|$ is not divisible by $4$. We then have that
 \[ \rank_{\La}\Big(e_{\eta}X^{\s}(E/\mF_{\infty})\Big) = \rank_{\La}\Big(e_{\bar{\eta}}X^{\s}(E/\mF_{\infty})^{\iota}\Big),\]
 and that the torsion $\La$-submodules of $e_{\eta}X^{\s}(E/\mF_{\infty})$ and $e_{\bar{\eta}}X^{\s}(E/\mF_{\infty})^{\iota}$ are pseudo-isomorphic.
\end{thm2}

Our result is in the spirit of prior results of Kim \cite{Kim08} (also see \cite{G89}). In fact, Kim has proven a variant of the above result for elliptic curves having good supersingular reduction all all primes above $p$. Our result thus is an improvement on his, namely, we allowed the presence of primes of good ordinary reduction above $p$ and our signed Selmer groups can have mixed signs.
As an application of our theorem, we prove the following.

\begin{cor2}[Corollary \ref{main corollary}]
Retain the setting of the above Theorem. Then we have $e_{\eta}X^{\s}(E/\mF_{\infty})=0$ if and only if $e_{\bar{\eta}}X^{\s}(E/\mF_{\infty})=0$.
\end{cor2}

The above theorem and corollary can be thought as a form of reflection principle reminiscent to that of class groups (for instance, see \cite[\S10.2]{Wa}) which relate an Iwasawa module in a given eigenspace
with another Iwasawa module in a ``reflected'' eigenspace. We also mention a viewpoint which may be interesting. Philosophically, one would expect that there should be a $p$-adic $L$-function on the analytic side which can be related to the signed Selmer groups via a main conjecture(s). Granted this, one will naturally hope for a result in the spirit of the Herbrand-Ribet Theorem which amounts to saying that the eigenspaces $e_{\eta}X^{\s}(E/\mF_{\infty})$ and $e_{\bar{\eta}}X^{\s}(E/\mF_{\infty})$ vanishes if and only if a certain twist of the conjectural $p$-adic $L$-function is a ``unit'' in an appropriate sense (see \cite[Corollary 5.2.7]{BCGKPST} or \cite{Wa}). For an elliptic curve $E$ defined over $\Q$, such a result is certainly true in view of recent progresses on the main conjecture (see \cite{K, Kob, SU, Wan}). However, a construction of $p$-adic $L$-function for more general base field is still very much absent. Hence our algebraic result can be viewed as a partial support to the above philosophy. It would definitely be interesting to look for Herbrand-Ribet Theorem type phenomenon in other Iwasawa-theoretical situations. We hope to pursue this line of study in a future work.

We now give an outline of our paper. In Section \ref{Algebraic preliminaries}, we collect various results concerning $\La$-modules and twist of $\La[G]$-modules for a cyclic group $G$ whose order is coprime to $p$. In Section \ref{local calculations}, we collect various facts on the arithmetic of an elliptic curve over a local field which will be needed for the subsequent discussion of the paper. Section \ref{Selmer} is where we introduce the signed Selmer groups. We also need to work with the so-called strict signed Selmer group which coincides with the signed Selmer group on the infinite level. Note that the strict signed Selmer group and the signed Selmer group need not coincide on the intermediate subextensions. Despite this, they enjoy nice duality properties which will be crucial in proving our results. Once the machinery of the strict signed Selmer groups is set up, we can give the proof of Theorem \ref{main theorem} and Corollary \ref{main corollary}.

\begin{ack}
The research of this article took place when S. Ahmed was a postdoctoral fellow at Central China Normal University, and he would like to acknowledge the hospitality
and conducive working conditions provided by the said institute. M. F. Lim gratefully acknowledges support by the
National Natural Science Foundation of China under Grant No. 11550110172 and Grant No. 11771164.
 \end{ack}

\section{Algebraic preliminaries} \label{Algebraic preliminaries}

\subsection{$\La$-modules}

Throughout the paper, $\La$ will denote the ring $\Zp\ps{T}$. By the Weierstrass preparation theorem (cf. \cite[(5.3.4)]{NSW} or  \cite[Theorem 7.3]{Wa}), every $f\in \La$ can be written as a product $u\cdot f_1$, where $u$ is a unit in $\La$ and $f_1$ is a Weierstrass polynomial in $\Zp[T]$, where one recalls that a polynomial $T^n +
c_{n-1}T^{n-1} + \cdots +c_0$ in $\Zp[T]$ is said to be a
Weierstrass polynomial if $p$ divides $c_i$ for every
$0\leq i \leq n-1$. We shall then define the degree of $f$ to be the degree of the polynomial $f_1$ and denote it by $\deg(f)$.
For a finitely generated $\La$-module $M$, there exist irreducible Weierstrass polynomials $f_j$, numbers $r$, $\al_i$, $\be_j$ and a homomorphism
\[ M \lra \La^r \oplus \Big(\bigoplus_{i=1}^s\La/p^{\al_i}\Big) \oplus \Big(\bigoplus_{j=1}^t\La/f_j^{\be_j}\Big) \]
with finite kernel and cokernel, where the numbers $r$, $\al_i$ $\be_j$ and the Weierstrass polynomials $f_j$ are determined by $M$ (cf. \cite[(5.3.8)]{NSW} or \cite[Theorem 13.12]{Wa}).
Then the rank of $M$ is given by the number $r$ which we denote by $\rank_{\La}(M)$. The $\mu$-invariant $\mu(M)$ (resp., $\la$-invariant $\la(M)$) of $M$ is defined to be $\sum_{i=1}^s\al_i$
(resp., $\sum_{j=1}^t\be_j\deg(f_j)$). Finally, the characteristic ideal of $M$ is defined to be the ideal generated by $\displaystyle\prod_{i=1}^sp^{\al_i}\prod_{j=1}^tf_j^{\be_j}$. We shall write $\Char(M)$ for this ideal. For a finitely generated $\La$-module, we write $M_{\La-tor}$ for its torsion $\La$-submodule. It is clear from the definition that one has $\Char(M)= \Char(M_{\La-tor})$.

\bl \label{rank compare}
Let $M$ and $N$ be two finitely generated $\La$-modules. Suppose that for each $m$,
the quantity $\big|(M/p^m)_{\Ga_n}\big|\Big/\big|(N/p^m)_{\Ga_n}\big|$ is bounded independent of $n$.
Then we have that $\rank_{\La}(M) = \rank_{\La}(N)$, and that $M[p^{\infty}]$ and $N[p^{\infty}]$ are pseudo-isomorphic as $\La$-modules.
\el

\bpf
For each fixed $m$, it follows from \cite[Theorem 2.5.1]{LimComp} that
\[ \big|(M/p^m)_{\Ga_n}\big| = p^{\mu(M/p^m)p^n+ O(1)}\]
as $n$ varies. Therefore, it follows from this and the hypothesis of the lemma that $\mu(M/p^m) = \mu(N/p^m)$ for every $m$. The required conclusion now follows from an application of \cite[Proposition 2.4.7]{LimComp}.
 \epf

For a finitely generated $\La$-module $M$, one can endow it with a natural compact topology (see \cite[Proposition 5.2.17]{NSW}). We then write $M^{\vee} = \Hom_{\cts}(M,\Qp/\Zp)$ for the group of continuous group homomorphisms from $M$ to $\Qp/\Zp$, where $\Qp/\Zp$ is given the discrete topology.

\bl \label{dist poly compare}
Let $M$ be a finitely generated $\La$-module. Suppose that \[ M \lra \La^r \oplus \Big(\bigoplus_{i=1}^s\La/p^{\al_i}\Big) \oplus \Big(\bigoplus_{j=1}^t\La/f_j^{\be_j}\Big) \]
with finite kernel and cokernel, where the numbers $r$, $\al_i$ $\be_j$ and the Weierstrass polynomials $f_j$ is determined by $M$. Then for every irreducible Weierstrass polynomial $f$, we have
\[ \corank_{\Zp}\big((M^{\vee}\ot \La/f^n)^\Ga\big) = \left(r + \sum_{j=1}^t\delta_{f_j, f}\min\{\be_j,n\}\right)\deg(f),\]
where $\delta_{f_j, f}=1$ or $0$ according as $f_j\La = f\La$ or $f_j\La\neq f\La$.
\el

\bpf
From the discussion in \cite[pp 109]{G89}, we have $(M^{\vee}\ot \La/f^n)^\Ga = \Hom_{\La}(M, \Qp/\Zp\ot_{\Zp}\La/f^n)$, where the $\Zp$-corank of the latter is the $\Zp$-rank of $\Hom_{\La}(M, \La/f^n)$. Since $\Hom_{\La}(-, \La/f^n)$ preserves finite direct sums and sends finite modules to finite modules, it suffices to prove the asserted equality in the lemma for an elementary $\La$-module which is a straightforward verification.
\epf

We can now state the following proposition which is essentially taken from \cite[Section 3]{G89}.

\bp \label{compare La modules}
Let $M$ and $N$ be two finitely generated $\La$-modules. Suppose that the following two statements hold.
\begin{enumerate}
\item[$(1)$] For each $m$,
the quantity $\big|(M/p^m)_{\Ga_n}\big|\Big/\big|(N/p^m)_{\Ga_n}\big|$ is bounded independent of $n$.
\item[$(2)$] For each irreducible Weierstrass polynomial $f$, we have
\[ \corank_{\Zp}\big((M^{\vee}\ot \La/f^n)^\Ga\big) =  \corank_{\Zp}\big((N^{\vee}\ot \La/f^n)^\Ga\big)\] for $n\geq 0$.
\end{enumerate}
Then we have $\rank_{\La}(M) = \rank_{\La}(N)$, and $M_{\La-tor}$ is pseudo-isomorphic to $N_{\La-tor}$.
\ep

\bpf
Lemma \ref{rank compare} tells us that $\rank_{\La}(M)= \rank_{\La}(N)$ and $M[p^{\infty}]$ is pseudo-isomorphic to $N[p^{\infty}]$. Hence, it remains to show that
for each irreducible Weierstrass polynomial $f$, the power of $f$ occurring in $M$ and $N$ are the same, which is precisely a consequence of Lemma \ref{dist poly compare}.
 \epf

\subsection{Twist of $\La$-modules}
Suppose that $G$ is cyclic of order coprime to $p$. (In most of our subsequent discussion, $G$ is a subgroup of $\Z/(p-1)\Z$.) Let $\eta$ be a fixed character of $G$. We write $e_{\eta}$ for the idempotent $\frac{1}{|G|}\sum_{\si\in G} \eta(\si)\si^{-1}\in\Zp[G]$. As before, write $\La$ for the power series ring $\Zp\ps{T}$ in one variable. We shall also view $e_{\eta}$ as an element of $\La[G]$. For a $\La[G]$-module $M$, we write $M^{\eta} = e_{\eta}M$. It is an easy exercise to check that $\si\in G$ acts on $M^{\eta}$ by multiplication by $\eta(\si)$.

We shall write $\Zp(\eta)$ for the $\Zp[G]$-module which is $\Zp$ as $\Zp$-module with the action of $G$ given by $\si\cdot x = \eta(\si)x$ for $\si\in G, x\in \Zp$. If $M$ is a $\La[G]$-module, we write $M(\eta) = M\ot_{\Zp}\Zp(\eta)$. (Warning: $M^{\eta}\neq M(\eta)$ in general; see Lemma \ref{twist and eta}(4) below.) For an element $f$ of $\La$ which is not divisible by $p$, we write $M_f = M\ot_{\Zp}\La/f$. We shall also denote by $f^{\iota}$ the power series $\displaystyle f\left(\frac{1}{1+T}-1\right)$. For ease of discussion, we record a simple lemma.

\bl \label{submodule sum}
Let $I$ be a finite indexing set. Suppose that $M_i$ and $N_i$ are $\Zp$-submodules of $M$. Assume that $N_i\subseteq M_i$ for every $i\in I$ and that
\[ M = \bigoplus_{i}M_i = \bigoplus_{i}N_i.\]
Then $M_i = N_i$ for every $i$.
\el

\bpf
 Let $x \in M_j$. Then viewing $x$ as an element of $M = \bigoplus_{i}N_i$, we have $ x = \sum_i y_i$ for $y_i\in N_i$. But since $N_i\subseteq M_i$ for every $i\in I$, the sum can be viewed as a direct sum decomposition in $M = \bigoplus_{i}M_i$. Hence we have $x =y_j$ which implies that $x\in N_j$. This proves the lemma.
\epf

We can now state the following lemma which tells us how the twists by $\eta$ and $f$ interact.

\bl \label{twist and eta}
Let $M$ be a $\La[G]$-module. Let $\eta$ be a character of $G$ and $f$ a power series of $\La$ which is not divisible by $p$. Viewing $\La/f$ as a $\La[G]$-module with a trivial $G$-action, we have the following equalities of $\La$-modules.
\begin{enumerate}
 \item[$(1)$] $e_{\eta}(M_f) = (e_{\eta}M)_f$.
 \item[$(2)$] $e_{\eta}(M^{\vee})_f = \Big(e_{\bar{\eta}}M_{f^{\iota}}\Big)^{\vee}$.
 \item[$(3)$] $e_{\eta}M_{div} = \big(e_{\eta}M\big)_{div}$.
 \item[$(4)$] $e_{\eta}M = \big(M(\bar{\eta})\big)^G$.
\end{enumerate}
\el

\bpf
 We begin by proving statement (1). Note that
 \[ M_f = \bigoplus_{\eta} e_{\eta} (M_f),\]
 where $\eta$ runs through all the character of $G$. On the other hand, we have
 \[ M_f = \left(\bigoplus_{\eta} e_{\eta}M\right)_f = \bigoplus_{\eta} \left(e_{\eta}M\right)_f.\]
 It is easy to verify that $(e_{\eta}M)_f\subseteq e_{\eta} (M_f)$. Hence the equality now follows from Lemma \ref{submodule sum}. Statements (2) and (3) can be proven similarly. For statement (4), one observes that $x\in e_{\eta}M \Leftrightarrow \si x = \eta(\si)x \Leftrightarrow \bar{\eta}(\si)\si x = x  \Leftrightarrow x\in \big(M(\bar{\eta})\big)^G$.
\epf

\section{Local considerations} \label{local calculations}

\subsection{Local duality} \label{Local duality subsection}

Let $K$ be a local field of characteristic zero, and $L$ a finite cyclic Galois extension of $K$ with Galois group $G$. Suppose that the order of $G$ is coprime to $p$. (In most of our subsequent discussion, $G$ is a subgroup of $\Z/(p-1)\Z$.) Let $\eta$ be a character of $G$. In this subsection, $M$ will always denote a finite $\La[\Gal(\bar{K}/K)]$-module of order a power of $p$.

\bl \label{coprime H1}
We have an isomorphism $e_{\eta}H^1(L,M) \cong H^1(K,M(\bar{\eta}))$.
\el

\bpf
By Lemma \ref{twist and eta}(4), we have $e_{\eta}H^1(L,M) = \big(H^1(L,M)(\bar{\eta})\big)^G$. Since $\Gal(\bar{L}/L)$ acts trivially on $\Zp(\bar{\eta})$, the latter is $H^1\big(L,M(\bar{\eta})\big)^G$.
By the inflation-restriction sequence, we have
\[ 0\lra H^1\left(G, M(\bar{\eta})^{\Gal(\bar{L}/L)}\right)\lra H^1\big(K,M(\bar{\eta})\big) \lra H^1\big(L,M(\bar{\eta})\big)^G \lra H^2\left(G, M(\bar{\eta})^{\Gal(\bar{L}/L)}\right).\]
Since $G$ is assumed to have order coprime to $p$, the $H^i(G,-)$ terms vanish, thus yielding the required isomorphism of the lemma.
\epf

Recall that the local Tate duality gives us a perfect pairing
\[ H^1(L,M) \times H^1(L,M^{*}) \lra \Qp/\Zp, \]
where $M^{*}=\Hom_{\Zp}(M,\mu_{p^{\infty}})$.

\bl \label{orthogonal complement}
Let $C$ be a subgroup of $H^1(L,M)$, and let $D$ be its orthogonal complement with respect to the above perfect pairing. For a given character $\eta$ of $G$ and an element $f\in \La$ which is coprime to $p$, $e_{\eta}C_f$ (resp., $e_{\bar{\eta}}D_{f^{\iota}}$) can be viewed as subgroups of $H^1(K,M(\bar{\eta})_f)$ (resp. $H^1(K,M^*(\eta)_{f^{\iota}})$). Furthermore, they are orthogonal complement to each other with respect to the
perfect pairing
\[ H^1(K,M(\bar{\eta})_{f}) \times H^1(K,M^*(\eta)_{f^{\iota}}) \lra \Qp/\Zp. \]
\el

\bpf
By the orthogonal relation between $C$ and $D$, we have short exact sequences
\[ 0\lra C\lra H^1(L,M) \lra D^{\vee}\lra 0, \]
\[ 0\lra D\lra H^1(L,M^*) \lra C^{\vee}\lra 0. \]
Taking idempotent and tensoring with $\ot_{\Zp}\La/f$ (or $\ot_{\Zp}\La/f^{\iota}$), we obtain
\[ 0\lra e_{\eta}C_f\lra H^1(K,M(\bar{\eta})_f) \lra (e_{\bar{\eta}}D_{f^{\iota}})^{\vee}\lra 0, \]
\[ 0\lra e_{\bar{\eta}}D_{f^{\iota}}\lra H^1(K,M^*(\eta)_{f^{\iota}}) \lra (e_{\eta}C_f)^{\vee}\lra 0, \]
where the middle terms of the short exact sequences are obtained via Lemma \ref{coprime H1}. This proves the claims of the lemma.
\epf

\subsection{Ordinary reduction at $p$} \label{ordinary subsection}

Let $K$ be a finite extension of $\Qp$, and write $L=K(\mu_p)$. Denote by $G$ the Galois group $\Gal(L/K)$ which has order dividing $p-1$. Let $K_{\infty}$ be the cyclotomic $\Zp$-extension of $K$ with $K_n$ the intermediate subfield with $|K_n:K|=p^n$. We write $L_{\infty}=K_{\infty}L$ which is the cyclotomic $\Zp$-extension of $L$ with intermediate subfield $L_n$ such that $|L_n:L|=p^n$. We identify the completed group algebras $\Zp\ps{\Gal(L_{\infty}/L)} = \Zp\ps{\Gal(K_{\infty}/K)} = \La$ and identify the Galois groups $\Gal(L_n/K_n)$ with $G=\Gal(L/K)$ for every $n$.

Now let $E$ denote an elliptic curve over $K$ with good ordinary reduction. Write $\widehat{E}$ (resp. $\widetilde{E}$) for its formal group (resp. reduced curve) at $K$. By \cite[Corollary 3.2]{CG}, we have that $H^1(L_{\infty}, \widehat{E}[p^{\infty}]) =0$ which in turn yields
\[ 0\lra \widehat{E}(L_{\infty})[p^{\infty}]\lra E(L_{\infty})[p^{\infty}] \lra \widetilde{E}(L_{\infty})[p^{\infty}]\lra 0. \]
The main theorem of Imai \cite{Im} says that $E(L_{\infty})[p^{\infty}]$ is finite. Consequently, so are $\widehat{E}(L_{\infty})[p^{\infty}]$ and $\widetilde{E}(L_{\infty})[p^{\infty}]$. In fact, we have the following.

\bl \label{Imai finite}
 Let $\eta$ be a character of $\Gal(L_n/K_n)$ and $f\in\La$ which is coprime to $p$. Then $\widehat{E}(K_{\infty})[p^{\infty}](\eta)_f$, $E(K_{\infty})[p^{\infty}](\eta)_f$ and $\widetilde{E}(K_{\infty})[p^{\infty}](\eta)_f$ are finite.
\el

\bpf
 Since $E(K_{\infty})[p^{\infty}](\eta)_f = E(K_{\infty})[p^{\infty}]\ot_{\Zp}\Zp(\eta)\ot_{\Zp}\La/f\La$, the finiteness follows from the above observation that $E(K_{\infty})[p^{\infty}]$ is finite. The remaining cases can be proven similarly.
\epf

Let $\e_n$ denote the natural map
\[ H^1\big(L_n, \widehat{E}[p^{\infty}]\big) \lra H^1\big(L_n, \Ep\big). \]
It follows from \cite[Proposition 4.5]{CG} that $(\im \e_n)_{div} = E(L_n)\ot\Qp/\Zp$.
Now define $C(E[p^m]/L_n)$ to be the kernel of the following map
\[ H^1\big(L_n, E[p^m]\big) \lra \frac{H^1\big(L_n, \Ep\big)}{(\im \e_n)_{div}}.\]
It is easy to see that $C(E[p^m]/L_n)$ identifies with $E(L_n)\ot\Z/p^m$. It then follows from this and a classical result of Tate \cite{Ta} that $C(E[p^m]/L_n)$ is orthogonal to itself under the
perfect pairing
\[ H^1(L_n,E[p^m]) \times H^1(L_n, E[p^m]) \lra \Qp/\Zp. \]

\bl
\label{ordinary orthogonal}
Fix a character $\eta$ of $\Gal(L_n/K_n)$ and $f\in\La$ which is coprime to $p$. Then $e_{\eta}C(E[p^m]/L_n)_f$ is orthogonal to $e_{\bar{\eta}}C(E[p^m]/L_n)_{f^{\iota}}$ under the
perfect pairing
\[ H^1\big(K_n, E[p^{m}](\bar{\eta})_f\big) \times
H^1\big(K_n, E[p^{m}](\eta)_{f^{\iota}}\big) \lra \Qp/\Zp. \]
\el

\bpf
This follows immediately from combining Lemma \ref{orthogonal complement} with the above discussion.
\epf

For later purposes, we require the following description of $e_{\eta}C(E[p^m]/L_n)_f$. Let $\theta_n$ denote the natural map
\[ H^1\big(K_n, \widehat{E}[p^{\infty}](\bar{\eta})_f\big) \lra H^1\big(K_n, \Ep(\bar{\eta})_f\big). \]
For each $m$, we define $C(E[p^m](\bar{\eta})_f/K_n)$ to be the kernel of the following map
\[ H^1\big(K_n, E[p^m](\bar{\eta})_f\big) \lra \frac{H^1\big(K_n, \Ep(\bar{\eta})_f\big)}{(\im \theta_n)_{div}}.\]

\bl \label{description of eC}
Let $\eta$ be a character of $\Gal(L_n/K_n)$ and $f\in\La$ which is coprime to $p$. Then we have \[e_{\eta}C(E[p^m]/L_n)_f\cong C(E[p^m](\bar{\eta})_f/K_n).\]
\el

\bpf
This follows from a straightforward verification by making use of Lemma \ref{twist and eta}.
\epf

The following lemma is the main calculation that will be needed in the eventual proof of our main theorem.

\bl \label{local ordinary calculation}
For a character $\eta$ of $\Gal(L_n/K_n)$ and $f\in\La$ which is coprime to $p$, we have
\[ \frac{\Big|C(E[p^m](\bar{\eta})_f/K_n)\Big|\Big|H^2\big(K_n, E[p^m](\bar{\eta})_f\big)\Big|}{\Big|H^1\big(K_n, E[p^m](\bar{\eta})_f\big)\Big|} =\frac{\big| E(K_n)[p^\infty](\bar{\eta})_f/p^m\big| \big|\widehat{E}(K_n)[p^{\infty}](\bar{\eta})_f/p^m\big|}{\big|E(K_n)[p^m](\bar{\eta})_f\big|}
p^{-m|K_n:\Qp|\deg(f)} .\]
\el

\bpf
  By the local Euler characteristics formula (cf. \cite[Theorem 7.3.1]{NSW}), we have
  \[ \frac{\Big|H^2\big(K_n, E[p^m](\bar{\eta})_f\big)\Big|}{\Big|H^1\big(K_n, E[p^m](\bar{\eta})_f\big)\Big|} = \frac{1}{\Big|H^0\big(K_n, E[p^m](\bar{\eta})_f\big)\Big|}\big|E[p^m](\bar{\eta})_f\big|^{-|K_n:\Qp|} = \frac{1}{\big|E(K_n)[p^m](\bar{\eta})_f\big|}p^{-2m|K_n:\Qp|\deg(f)}. \]
  On the other hand, by Lemma \ref{description of eC}, we have a commutative diagram
  \[   \xymatrixrowsep{0.35in}
\entrymodifiers={!! <0pt, .8ex>+} \SelectTips{eu}{}\xymatrix{
    0 \ar[r]^{} & C(E[p^m](\bar{\eta})_f/K_n) \ar[d] \ar[r] &  H^1\big(K_n, E[p^m](\bar{\eta})_f\big)
    \ar[d] \ar[r] & \displaystyle \frac{H^1\big(K_n, \Ep(\bar{\eta})_f\big)}{(\im \theta_n)_{div}}[p^m]\ar@{=}[d]  \ar[r] &0\\
    0 \ar[r]^{} & (\im \theta_n)_{div}[p^m] \ar[r]^{} & H^1\big(K_n, E[p^\infty](\bar{\eta})_f\big)[p^m] \ar[r] & \displaystyle\frac{H^1\big(K_n, \Ep(\bar{\eta})_f\big)}{(\im \theta_n)_{div}}[p^m] \ar[r] & 0 } \]
with exact rows, where the surjectivity of the rightmost arrow on the top row follows from the fact that the rightmost arrow on the bottom row and the middle vertical map are surjective. It then follows from the snake lemma that we have a short exact sequence
 \[ 0 \lra E(K_n)[p^\infty](\bar{\eta})_f/p^m  \lra  C(E[p^m](\bar{\eta})_f/K_n) \lra (\im \theta_n)_{div}[p^m] \lra 0 \]
 which in turn yields an equality
 \[ \big| C(E[p^m](\bar{\eta})_f/K_n) \big| = \big| E(K_n)[p^\infty](\bar{\eta})_f/p^m\big| \big| (\im \theta_n)_{div}[p^m] \big|. \]
 We now calculate $\big| (\im \theta_n)_{div}[p^m] \big|$. To do this, we observe that there is an exact sequence
\[ \widetilde{E}(K_n)[p^{\infty}](\bar{\eta})_f \lra H^1\big(K_n, \widehat{E}[p^{\infty}](\bar{\eta})_f\big) \stackrel{\theta_n}{\lra} H^1\big(K_n, \Ep(\bar{\eta})_f\big). \]
Since $\widetilde{E}(K_n)[p^{\infty}](\bar{\eta})_f$ is finite, it follows that
\[ \corank_{\Zp}\big( (\im \theta_n)_{div}\big) = \corank_{\Zp}\big(\im \theta_n \big)=
\corank_{\Zp} H^1\big(K_n, \widehat{E}[p^{\infty}](\bar{\eta})_f\big).\]
By \cite[pp 111, Formula (23)]{G89}, the latter quantity is precisely \[|K_n:\Qp|\corank_{\Zp}\widehat{E}[p^{\infty}](\bar{\eta})_f = |K_n:\Qp|\deg f.\] Hence we have $\big| (\im \theta_n)_{div}[p^m] \big| = p^{m|K_n:\Qp|\deg f}$. Combining all the calculations, we obtain the desired formula of the lemma.
\epf

\subsection{Supersingular reduction at $p$} \label{supersingular subsection}

In this subsection, $K$ will denote a finite unramified extension of $\Qp$, and $E$ is an elliptic curve defined over $\Qp$ which is assumed to have good supersingular reduction with $a_p=0$. Set $L_n =K(\mu_{p^{n+1}})$ and $L_{\infty}=\cup_n L_n$. Write $\Ga_n=\Gal(L_{\infty}/L_n)$ for $n\geq 0$. We may identify $\Gal(L/K)$ with a subgroup of $\Gal(L_n/K)$ for every $n =1,...,\infty$. We then write $K_n$ for the fixed field of this subgroup.

\bl \label{supersingular points}
 Retain settings as above. Then $\widehat{E}(L_{n})$ and $E(L_n)$ have no $p$-torsion for every $n\geq -1$.
 \el

\bpf
The first assertion follows from \cite[Proposition 3.1]{KO}. The second follows from the first, the following short exact sequence
\[ 0\lra \widehat{E}(L_n)\lra E(L_n) \lra \widetilde{E}(\mathbb{F}_{p^d})\lra 0\]
and that $\widetilde{E}(\mathbb{F}_{p^d})$ has no $p$-torsion by the supersingular assumption. \epf

\bl \label{supersingular points twist}
Let $\eta$ be a character of $\Gal(L_n/K_n)$ and $f$ a distinguished polynomial in $\La$.  Then $\widehat{E}(K_{n})[p^{\infty}](\eta)_f = 0$ for every $n\geq 0$.
 \el

\bpf
This is immediate from Lemma \ref{supersingular points}. \epf

Following \cite{Kim07, Kim08, KimPM, Kim14, KO, Kob}, we define
the following groups
\[\widehat{E}^+(L_{n}) = \{ P\in \widehat{E}(L_{n})~:~\mathrm{tr}_{n/m+1}(P)\in \widehat{E}(L_{m}), 2\mid m, 0\leq m \leq n-1\}, \]
\[\widehat{E}^-(L_{n}) = \{ P\in \widehat{E}(L_{n})~:~\mathrm{tr}_{n/m+1}(P)\in \widehat{E}(L_{m}), 2\nmid m, -1\leq m \leq n-1\}, \]
where $\mathrm{tr}_{n/m+1}: E(L_{n}) \lra E(L_{m+1})$ denotes the trace map and $L_{-1}$ is understood to be $K$.  We also write
\[ \delta =\begin{cases} 0,  & \mbox{if $|K:\Qp|\neq 0$ (mod 4) or $\eta\neq 1$}, \\
2, & \mbox{otherwise}.\end{cases}\]

\bp \label{Kita-O results}
Retain settings as above. The following statements are then valid.
\begin{enumerate}
\item[$(a)$] $\big(\widehat{E}^{\pm}(L_{\infty})^{\eta}\ot\Qp/\Zp\big)^{\Ga_n}$ are cofree $\Zp$-modules for all $n$, and
    \[ \corank_{\Zp}\left(\big(\widehat{E}^{s}(L_{\infty})^{\eta}\ot\Qp/\Zp\big)^{\Ga_n}\right)
    =\begin{cases} |K:\Qp|p^n+\delta,  & \mbox{if $s=+$}, \\
 |K:\Qp|p^n, & \mbox{if $s=-$}.\end{cases}\]
\item[$(b)$] One has \[\left(\widehat{E}^{s}(L_{\infty})^{\eta}\ot\Qp/\Zp\right)^{\vee}\cong
    \begin{cases} \Zp\ps{X}^{\oplus |K:\Qp|}\oplus \Zp^{\oplus \delta},  & \mbox{if $s=+$}, \\
 \Zp\ps{X}^{\oplus |K:\Qp|}, & \mbox{if $s=-$}. \end{cases}\]
\end{enumerate}
\ep

\bpf
Statement (a) is \cite[Corollary 3.25]{KO}, and statement (b) follows from \cite[Theorem 3.34]{KO} (also see \cite[Theorems 2.7 and 2.8]{Kim14}).
\epf

Write $\mathbb{H}^{\pm}_{\infty}= \hat{E}^{\pm}(L_{\infty})\ot\Qp/\Zp$ and $\mathbb{H}^{\pm}_n= \left(\hat{E}^{\pm}(L_{\infty})\ot\Qp/\Zp\right)^{\Ga_n}$.
By the Hochshild-Serre spectral sequence and Lemma \ref{supersingular points}, we have an isomorphism
\[ H^1(L_n,\Ep)\cong H^1(L_{\infty},\Ep)^{\Ga_n}.\]
Via this isomorphism, we may view $\mathbb{H}^{\pm}_n$ as a subgroup of $H^1(L_n,E(p))$.
Let $M_n^{\pm}$ be the exact annihilator of $\mathbb{H}^{\pm}_n$ with respect to the local Tate pairing
\[ H^1(L_n, \Ep)\times H^1(L_n, T_pE)\lra \Qp/\Zp. \]
In other words, we have short exact sequences
\[ 0\lra \mathbb{H}^{\pm}_n \lra H^1(L_n, \Ep) \lra \big(M_n^{\pm}\big)^{\vee}\lra 0, \]
\[ 0\lra M_n^{\pm} \lra H^1(L_n, T_pE) \lra \big(\mathbb{H}^{\pm}_n\big)^{\vee}\lra 0. \]
Since $\mathbb{H}^{\pm}_n$ is divisible by Proposition \ref{Kita-O results}(a), it follows that there are two short exact sequences
\[ 0\lra \mathbb{H}^{\pm}_n[p^m] \lra H^1(L_n, \Ep)[p^m] \lra \big(M_n^{\pm}/p^m\big)^{\vee}\lra 0, \]
\[ 0\lra M_n^{\pm}/p^m \lra H^1(L_n, T_pE)/p^m \lra \big(\mathbb{H}^{\pm}_n[p^m]\big)^{\vee}\lra 0. \]

By appealing to Lemma \ref{supersingular points} again, one has $H^1(L_n, \Ep)[p^m]\cong H^1(L_n, E[p^m])$ and $H^1(L_n, T_pE)/p^m \cong H^1(L_n, E[p^m])$. Combining this with the above short exact sequences, we see that $M_n^{\pm}/p^m$ is the exact annihilator of $\mathbb{H}^{\pm}_n[p^m]$ with respect to the local Tate pairing
\[ H^1(L_n, E[p^m])\times H^1(L_n, E[p^m])\lra \Z/p^m. \]

In fact, one even has the following.

\bp[Kim] \label{Kim duality}
One always has $\mathbb{H}^{-}_n[p^m] = M_n^{-}/p^m$. In other words, $\mathbb{H}^{-}_n[p^m]$ is the exact annihilator of itself with respect to the local Tate pairing
\[ H^1(L_n, E[p^m])\times H^1(L_n, E[p^m])\lra \Z/p^m. \]
If $|K:\Qp|$ is not divisible by $4$, we then have the same conclusion for $\mathbb{H}^{+}_n[p^m]$.
\ep

\bpf
The first assertion is established by Kim in \cite[Proposition 3.15]{Kim07}. One can check that the same proof carries over for the second assertion under the assumption that $|K:\Qp|$ is not divisible by $4$ (also see \cite[Proposition 3.4]{KimPM} or \cite[Theorem 2.9]{Kim14}).
\epf

\bc \label{supersingular orthogonal}
Write $K_n = L_n^G$. Then
$\mathbb{H}^{-}_n[p^m]_f^{\eta}$ is the orthogonal complement to $\mathbb{H}^{-}_n[p^m]_{f^{\iota}}^{\bar{\eta}}$ with respect to the perfect pairing
\[ H^1(K_n, E[p^m](\bar{\eta})_f)\times H^1(K_n, E[p^m](\eta)_{f^{\iota}})\lra \Z/p^m. \]
If $|K:\Qp|$ is not divisible by $4$, we then have the same conclusion with $\mathbb{H}^{+}_n[p^m]_f^{\eta}$ and $\mathbb{H}^+_n[p^m]_{f^{\iota}}^{\bar{\eta}}$.
\ec

\bpf
This follows from a combination of Lemma \ref{orthogonal complement} and Proposition \ref{Kim duality}.
\epf

We end the subsection with a preliminary calculation which will be required in the proof of our main theorem.

\bl \label{local supersingular lemma}
Writing $K_n = L_n^G$, we have
\[ \frac{\big|\mathbb{H}^{-}_n[p^m]_f^{\eta}\big|~\big|H^2(K_n, E[p^m](\bar{\eta})_f\big|}{\big|H^1(K_n, E[p^m](\bar{\eta})_f\big|} = p^{-m|K_n:\Qp|\deg(f)}. \]
If $|K:\Qp|$ is not divisible by $4$, we have the same conclusion with $\mathbb{H}^{+}_n[p^m]_f^{\eta}$.
\el

\bpf
 By the local Euler characteristics formula (cf. \cite[Theorem 7.3.1]{NSW}), we have
  \[ \frac{\Big|H^2\big(K_n, E[p^m](\bar{\eta})_f\big)\Big|}{\Big|H^1\big(K_n, E[p^m](\bar{\eta})_f\big)\Big|} = \frac{1}{\Big|H^0\big(K_n, E[p^m](\bar{\eta})_f\big)\Big|}\big|E[p^m](\bar{\eta})_f\big|^{-|K_n:\Qp|} = p^{-2m|K_n:\Qp|\deg(f)}, \]
  where the final equality follows from Lemma \ref{supersingular points twist}. Now observe that
  \[ \mathbb{H}^{\pm}_n[p^m]_f^{\eta} = \left(e_{\eta}\left(\hat{E}^{\pm}(L_{\infty})\ot\Qp/\Zp\right)^{\Ga_n}\right)_f  = \left(\left(\hat{E}^{\pm}(L_{\infty})^{\eta}\ot\Qp/\Zp\right)^{\Ga_n}\right)_f\]
  By Proposition \ref{Kita-O results}(a) (and assuming that $|K:\Qp|$ is not divisible by $4$ in the $+$ case), the latter has size $p^{m|K:\Qp|\deg(f)}$. The conclusion of the lemma now follows from these calculations.
  \epf

\subsection{Away from $p$} \label{away from p subsection}

In this subsection, $l$ is a prime $\neq p$. Let $K$ be a finite extension of $\Q_l$, and write $L=K(\mu_p)$. Denote by $G$ the Galois group $\Gal(L/K)$ which has order dividing $p-1$. Let $K_{\infty}$ be the cyclotomic $\Zp$-extension of $K$. We write $L_{\infty}=K_{\infty}L$ which is the cyclotomic $\Zp$-extension of $L$. Similarly as before, we identify the completed group algebras $\Zp\ps{\Gal(L_{\infty}/L)} = \Zp\ps{\Gal(K_{\infty}/K)} = \La$ and identify the Galois groups $\Gal(L_n/K_n)$ with $\Gal(L/K)$ for every $n$.

For any discrete $\Gal(\bar{K}/K)$-module $M$ and algebraic extension $\mathcal{K}$ of $K$, we write $H^1_{ur}(\mathcal{K}, M)$ for the kernel of the map
\[ H^1(\mathcal{K}, M) \lra H^1(\mathcal{K}^{ur}, M),\]
where $\mathcal{K}^{ur}$ is the maximal unramified extension of $\mathcal{K}$. Via the inflation-restriction sequence, $H^1_{ur}(\mathcal{K}, M)$ identifies with $H^1(\Gal(\mathcal{K}^{ur}/\mathcal{K}), M(\mathcal{K}^{ur}))$.

Now let $E$ denote an elliptic curve over $K$. Let $\eta$ be a character of $\Gal(L_n/K_n)$ and $f$ a distinguished polynomial in $\La$.

\bl \label{local away from p lemma}
We have $e_{\eta}H^1_{ur}(L_n, E[p^m])_f = H^1_{ur}(K_n, E[p^m](\bar{\eta})_f)$ and
\[ \frac{\big|H^1_{ur}(K_n, E[p^m](\bar{\eta})_f)\big|~\big|H^2(K_n, E[p^m](\bar{\eta})_f\big|}{\big|H^1(K_n, E[p^m](\bar{\eta})_f\big|} = 1. \]
\el

\bpf
The first equality follows from a direct calculation appealing to Lemmas \ref{twist and eta} and \ref{coprime H1}. Local Euler characteristics formula (cf. \cite[Theorem 7.3.1]{NSW}) then tells us that
  \[ \frac{\big|H^0(K_n, E[p^m](\bar{\eta})_f)\big|~\big|H^2(K_n, E[p^m](\bar{\eta})_f\big|}{\big|H^1(K_n, E[p^m](\bar{\eta})_f\big|} = 1. \]
  Finally, one checks easily by definition that $H^0(K_n, E[p^m](\bar{\eta})_f) = H^1_{ur}(K_n, E[p^m](\bar{\eta})_f)$. This proves the second equality.
\epf

\bl \label{orthogonal away from p lemma}
Fix a character $\eta$ of $\Gal(L_n/K_n)$ and $f\in\La$ which is coprime to $p$. Then $H^1_{ur}(K_n, E[p^m](\bar{\eta})_f)$ is orthogonal to $H^1_{ur}(K_n, E[p^m](\eta)_{f^{\iota}})$ under the
perfect pairing
\[ H^1\big(K_n, E[p^{m}](\bar{\eta})_f\big) \times
H^1\big(K_n, E[p^{m}](\eta)_{f^{\iota}}\big) \lra \Qp/\Zp. \]
\el

\bpf
It is well-known that $H^1_{ur}(L_n, E[p^m])$ is orthogonal to itself under the
perfect pairing
\[ H^1\big(L_n, E[p^{m}]\big) \times
H^1\big(L_n, E[p^{m}]\big) \lra \Qp/\Zp \]
(for instance, see \cite[pp 113]{G89}). The conclusion of the lemma now follows from Lemmas \ref{orthogonal complement} and \ref{local away from p lemma}.
\epf

\section{Mixed signed Selmer groups} \label{Selmer}

\subsection{Main result}
 As a start, we recall the definition of the signed Selmer group \cite{Kim07, Kim08, KimPM, Kim14, KO, Kob}. Let $F'$ be a number field and $E$ an elliptic curve defined over $F'$. Let $F$ be a finite extension of $F'$. Denote by $\Si_p$ for the set of primes of $F$ lying above $p$. We shall also write $\Si_p^{ord}$ (resp., $\Si_p^{ss}$) for the set of primes in $\Si_p$ at which $E$ has good ordinary reduction (resp., good supersingular reduction). The following assumptions will be in full force throughout the paper.

\begin{itemize}
\item[(S1)] The elliptic curve $E$ has good reduction at all primes of $F'$ above $p$, and at least one of which is a supersingular reduction prime of $E$.

 \item[(S2)] For each $u$ of $F'$ above $p$ at which the elliptic curve $E$ has supersingular reduction, assume that the following statements are valid.

 (a) $F'_u=\Qp$.

 (b) $a_u = 1 + p - |\tilde{E}_u(\mathbb{F}_p)| = 0$, where $\tilde{E}_u$ is the reduction of $E$ at $u$.

 (c) $u$ is unramified in $F/F'$.
\end{itemize}

Denote by $\mF_n = F(\mu_{p^{n+1}})$ and $\mF_{\infty} = \cup_n \mF_n$. By (S1)-(S2), $\Gal(\mF_0/F)\cong \Z/(p-1)\Z$ and every supersingular prime of $E$ above $p$ is totally ramified in $\mF_{\infty}/F$. Thus, we have $\Gal(\mF_{\infty}/F)\cong \Gal(\mF_{\infty}/\mF_0)\times \Gal(\mF_0/F) \cong \Zp\times \Z/(p-1)\Z$. We shall write $\Ga_n = \Gal(\mF_{\infty}/\mF_n)$ and $G= \Gal(\mF_{0}/F)$. We shall also write $\Ga= \Ga_0$. By fixing a topological generator of $\Ga$, we identify $\Zp\ps{\Gal(\mF_{\infty}/F)}$ with $\La[G]$.

It follows from (S1) and (S2) that for each prime $w$ of $F$ above $p$ at which $E$ has supersingular reduction, there is a unique prime of $\mF_n$ lying above the said prime which, by abuse of notation, is still denoted by $w$.
For such a $w$, we choose a sign $s_w$ for it and write $\s = (s_w)_{w\in \Si^{ss}_p}\in \{\pm\}^{\Si^{ss}_p}$.
 For each $s_w$, define $E^{s_w}(\mF_{n,w})$ as in Subsection \ref{supersingular subsection}.
The mixed signed Selmer group is then given by
\begin{align*}
    \Sel^{\overrightarrow{s}}(E/\mF_n):=\ker\Bigg(H^1(\mF_n,\Ep)
    \lra&\\ \bigoplus_{w\in\Si^{ss}(F_n)}\frac{H^1(\mF_{n,w},\Ep)}{\widehat{E}^{s_w}(\mF_{n,w})\ot\Qp/\Zp}
    \times &
    \bigoplus_{w\in\Si^{ord}(\mF_n)}
    \frac{H^1(\mF_{n,w},\Ep)}{E(\mF_{n,w})\ot\Qp/\Zp}\times \bigoplus_{w\nmid p} H^1(\mF_{n,w},\Ep)  \Bigg),
\end{align*}
Define $\Sel^{\s}(E/\mF_{\infty})=\ilim_n \Sel^{\overrightarrow{s}}(E/\mF_n)$.
Denote by $X^{\s}(E/\mF_{\infty})$ the Pontryagin dual of $\Sel^{\overrightarrow{s}}(E/\mF_\infty)$. It is not difficult to verify that this is finitely generated over $\Zp\ps{\Ga}$. In fact, one expects the following conjecture which is a natural extension of Mazur \cite{Maz} and Kobayashi \cite{Kob}.

\medskip \noindent \textbf{Conjecture.} $X^{\overrightarrow{s}}(E/\mF_{\infty})$ is a torsion $\Zp\ps{\Ga}$-module.

\medskip
When $E$ has good ordinary reduction at all primes above $p$, the above conjecture is precisely Mazur's conjecture \cite{Maz} which is known in the case when $E$ is defined over $\Q$ and $F$ an abelian extension of $\Q$ (see \cite{K}). For an elliptic curve over $\Q$ with good supersingular reduction at $p$, this conjecture was established by Kobayashi
(cf. \cite{Kob}; also see \cite{BL} for some recent progress on this conjecture).

From now on, let $\Sigma$ denote a fixed finite set of primes of $F$ containing those above $p$, the ramified primes of $F/F'$ and all the bad reduction primes of $E$. Write $F_{\Sigma}$ for the maximal algebraic extension of $F$ which is unramified outside $\Sigma$. For any (possibly infinite) extension $F\subseteq L\subseteq F_{\Sigma}$, write $G_{\Sigma}(L)= \Gal(F_{\Sigma}/L)$.
The signed Selmer group of $E$ over $\mF_\infty$ can then be equivalently defined by
\begin{align*}
    \Sel^{\overrightarrow{s}}(E/\mF_\infty)=\ker\Bigg(H^1(G_{\Sigma}(\mF_\infty),\Ep)
    \lra&\\ \bigoplus_{w\in\Si^{ss}(\mF_\infty)}\frac{H^1(\mF_{\infty,w},\Ep)}{E^{s_w}(\mF_{\infty,w})\ot\Qp/\Zp}\times
        &\bigoplus_{w\in\Si^{ord}(\mF_\infty)}\frac{H^1(\mF_{\infty,w},\Ep)}{E(\mF_{\infty,w})\ot\Qp/\Zp}
    \times\bigoplus_{w\in\Sigma, w\nmid p} H^1(\mF_{\infty,w},\Ep)  \Bigg),
\end{align*}

We can now state our main results.

\bt \label{main theorem}
 Let $\eta$ be a character of $G=\Gal(F(\mu_p)/F)$. For each $v\in \Si^{ss}$ with $s_v=+$, assume further that $|F_v:\Qp|$ is not divisible by $4$. We then have that
 \[ \rank_{\La}\Big(e_{\eta}X^{\s}(E/\mF_{\infty})\Big) = \rank_{\La}\Big(e_{\bar{\eta}}X^{\s}(E/\mF_{\infty})^{\iota}\Big),\]
 and that the torsion $\La$-submodules of $e_{\eta}X^{\s}(E/\mF_{\infty})$ and $e_{\bar{\eta}}X^{\s}(E/\mF_{\infty})^{\iota}$ are pseudo-isomorphic.
\et

\bc \label{main corollary}
Retain the setting of Theorem \ref{main theorem}. Then the following statements are valid.
\begin{enumerate}
\item[$(a)$] $e_{\eta}X^{\s}(E/\mF_{\infty})$ is torsion over $\La$ if and only if $e_{\bar{\eta}}X^{\s}(E/\mF_{\infty})$ is torsion over $\La$.
\item[$(b)$] $e_{\eta}X^{\s}(E/\mF_{\infty})=0$ if and only $e_{\bar{\eta}}X^{\s}(E/\mF_{\infty})=0$.
\end{enumerate}
\ec

The proof of the theorem will be given in Subsection \ref{proof of main theorem subsec} and the proof of the corollary will be given in Subsection \ref{proof of main corollary subsec}.

\subsection{Strict signed Selmer groups}

To prove our main result, we need to work with the so-called strict signed Selmer groups. The strict signed Selmer group of $E$ over $\mF_n$ is defined by
\[
    \Sel^{\overrightarrow{s},str}(E/\mF_n):=\ker\Bigg(H^1(G_{\Sigma}(\mF_n),\Ep)
    \lra \bigoplus_{w\in\Si^{ss}(\mF_n)}\frac{H^1(\mF_{n,w},\Ep)}
    {C_w\Big(\Ep/\mF_n\Big)}\Bigg),\]
where \[ C_w\Big(\Ep/\mF_n\Big) = \begin{cases}  \Big(E^{s_w}(\mF_{\infty,w})\ot\Qp/\Zp\Big)^{\Ga_n},& \mbox{if $w\in \Si_p(\mF_n)$}, \\
 \Bigg(\im\left(H^1(\mF_{n,w},\widehat{E}[p^{\infty}])\lra H^1(\mF_{n,w}, \Ep)\right)\Bigg)_{div},& \mbox{if $w\in \Si_p^{ord}(\mF_n)$}, \\
 H^1_{ur}(\mF_{n,w},\Ep). & \mbox{if $w\nmid p$,} \end{cases}\]
Define $\Sel^{\overrightarrow{s},str}(E/\mF_\infty):=\ilim_n\Sel^{\overrightarrow{s},str}(E/\mF_n)$.

Note that in the case when $\Si_p^{ss}=\emptyset$, this is the strict Selmer group as defined in the sense of Greenberg \cite{G89} (also see \cite{Guo}). Our choice of naming the above as the strict signed Selmer group is inspired by this observation.
In general, this strict Selmer group needs not agree with the signed Selmer group on the intermediate level $\mF_n$. But they do coincide upon taking limit which is the content of the next proposition.

\bp \label{strict = usual}
 We have an identification $\Sel^{\overrightarrow{s},str}(E/\mF_\infty) \cong \Sel^{\overrightarrow{s}}(E/\mF_\infty)$.
\ep

\bpf
 By \cite[Proposition 4.5]{CG}, we have $C_w(\Ep/\mF_n) = E(\mF_{n,w})\ot\Qp/\Zp$ for all $w\in \Si_p^{ord}(\mF_n)$ and so the local condition at $\Si_p^{ord}(\mF_n)$ for the strict signed Selmer group and the signed Selmer group at each $\mF_n$ (and hence $\mF_{\infty}$)) are the same. Similarly, one has the same conclusion for primes outside $p$. Let $w$ be a prime $\Si_p^{ss}$. Recall that we also write $w$ for the prime of $\mF_{\infty}$ above $w$. In general,
$\Big(E^{s_w}(\mF_{\infty,w})\ot\Qp/\Zp\Big)^{\Ga_n}$ and $E^{s_w}(\mF_{n,w})\ot\Qp/\Zp$ need not agree (one can see this by comparing their $\Zp$-coranks). But upon taking limit, they are both equal to
$E^{s_w}(\mF_{\infty,w})\ot\Qp/\Zp$. In conclusion, we have established our assertion.
\epf

We now show that this strict signed Selmer group enjoys a control theorem parallel to that in the ordinary situation (see \cite{G99, Maz}; also see \cite{Kim14, Kob}).

\bp \label{control theorem for strict}
There is an injection
\[ \Sel^{\s,str}(E/\mF_n) \lra \Sel^{\overrightarrow{s}}(E/\mF_\infty)^{\Ga_n}\]
with finite cokernel which is bounded independent of $n$.
\ep

\bpf
Consider the following diagram
\[   \xymatrixrowsep{0.25in}
\entrymodifiers={!! <0pt, .8ex>+} \SelectTips{eu}{}\xymatrix{
    0 \ar[r]^{} &\Sel^{\s,st}(E/\mF_n)  \ar[d] \ar[r] &  H^1\big(G_{\Sigma}(\mF_n),\Ep\big)
    \ar[d] \ar[r] & \displaystyle \bigoplus_{w\in\Si(\mF_n)}\frac{H^1\big(\mF_{n,w},\Ep\big)}
    {C_w\Big(\Ep/\mF_n\Big)}\ar[d]  \\
    0 \ar[r]^{} & \Sel^{\s,st}(E/\mF_\infty)^{\Ga_n} \ar[r]^{} & H^1\big(G_{\Sigma}(\mF_\infty),\Ep\big)^{\Ga_n} \ar[r] &
    \displaystyle\left(\bigoplus_{w\in\Si(F_\infty)}\frac{H^1\big(\mF_{\infty,w},\Ep\big)}
    {C_w\Big(\Ep/\mF_\infty\Big)}\right)^{\Ga_n}  } \]
with exact rows. Now the middle vertical map is surjective with kernel $H^1\big(\Ga_n, E(\mF_{\infty})[p^{\infty}]\big)$. By Lemma \ref{supersingular points} and (S1), this kernel is trivial. Hence it remains to show that the rightmost vertical map has finite kernel which is bounded independent of $n$. For primes in $\Si_p^{ord}$ and primes not dividing $p$, this is established in the mist of proving the control theorem in the ordinary case (for instance, see \cite[Theorem 1.2 and Section 3]{G99}). It therefore remains to consider the primes in $\Si_p^{ss}$. By our definition of $C_w\Big(\Ep/\mF_n\Big)$, we have
\[C_w\Big(\Ep/\mF_n\Big) =   \Big(E^{s_w}(\mF_{\infty,w})\ot\Qp/\Zp\Big)^{\Ga_n}  = C_w\Big(\Ep/\mF_\infty\Big)^{\Ga_n}\]
and hence the map is injective for these primes.
\epf

We now define mod-$p^m$ strict signed Selmer group on the level of $F_n$. This is given by
\[
    \Sel^{\overrightarrow{s},str}(E[p^m]/\mF_n):=\ker\Bigg(H^1(G_{\Sigma}(\mF_n),E[p^m])
    \lra \bigoplus_{w\in\Si^{ss}(\mF_n)}\frac{H^1(\mF_{n,w},E[p^m])}
    {C_w\Big(E[p^m]/\mF_n\Big)}\Bigg),\]
where \[ C_w\Big(E[p^m]/\mF_n\Big) = \begin{cases}  \Big(E^{s_w}(\mF_{\infty,w})\ot\Qp/\Zp\Big)^{\Ga_n}[p^m],& \mbox{if $w\in \Si_p^{ss}(\mF_n)$}, \\
 \ker\left(H^1(\mF_{n,w},E[p^{m}])\lra \displaystyle\frac{H^1(\mF_{n,w}, E[p^{\infty}])}{C_w\big(\Ep/\mF_{n}\big)}\right),& \mbox{if $w\in \Si_p^{ord}(\mF_n)$}, \\
 H^1_{ur}(\mF_{n,w},E[p^m]). & \mbox{if $w\nmid p$,} \end{cases}\]

\bp \label{control theorem for strict torsion}
The kernel and cokernel of the map
\[ \Sel^{\s,str}(E[p^m]/\mF_n) \lra \Sel^{\overrightarrow{s},str}(E/\mF_n)[p^m]\]
are finite and bounded independent of $m$ and $n$.
\ep

\bpf
Consider the following diagram
\[   \xymatrixrowsep{0.25in}
\entrymodifiers={!! <0pt, .8ex>+} \SelectTips{eu}{}\xymatrix{
    0 \ar[r]^{} &\Sel^{\s,str}(E[p^m]/\mF_n)  \ar[d] \ar[r] &  H^1\big(G_{\Sigma}(\mF_n),E[p^m]\big)
    \ar[d] \ar[r] & \displaystyle \bigoplus_{w\in\Si(\mF_n)}\frac{H^1\big(\mF_{n,w},E[p^m]\big)}
    {C_w\Big(E[p^m]/\mF_n\Big)}\ar[d]  \\
    0 \ar[r]^{} &  \Sel^{\overrightarrow{s},str}(E/\mF_n)[p^m] \ar[r]^{} & H^1\big(G_{\Sigma}(\mF_n),\Ep\big)[p^m] \ar[r] &
    \displaystyle\left(\bigoplus_{w\in\Si(\mF_n)}\frac{H^1\big(\mF_{n,w},\Ep\big)}
    {C_w\Big(\Ep/\mF_\infty\Big)}\right)[p^m]  } \]
with exact rows. By appealing to Lemma \ref{supersingular points}, one can show that the middle map is an isomorphism. It remains to show that the rightmost map is finite and bounded independent of $m$ and $n$. For primes at $\Si_p^{ord}(\mF_n)$, the map is injective by our definition of $C_w\Big(E[p^m]/\mF_n\Big)$. For primes at $\Si_p^{ss}(\mF_n)$ and primes not dividing $p$, the finiteness and boundness are shown in the proof of \cite[Lemma 3.4]{Kim08}. This finishes the proof of the proposition.
\epf

\subsection{Proof of the Theorem \ref{main theorem}} \label{proof of main theorem subsec}
This subsection is devoted to proving Theorem \ref{main theorem}. As a start, we determine $e_{\eta}\Sel^{\overrightarrow{s}}(E/\mF_\infty)_f$ for a character $\eta$ of $G=\Gal(F(\mu_p)/F)$ and a power series $f\in \La$ which is not divisible by $p$. For the remainder of the paper, we write $F_n = \mF_n^G$ and $F_{\infty} =\cup_nF_n$. Note that $F_{\infty}$ is precisely the cyclotomic $\Zp$-extension of $F$.

\bl \label{twist selmer}
For a character $\eta$ of $G=\Gal(F(\mu_p)/F)$ and a power series $f\in \La$ which is not divisible by $p$, we have
\[ e_{\eta}\Sel^{\overrightarrow{s},str}(E/\mF_n)_f\cong \ker\Bigg(H^1\big(G_{\Sigma}(F_n),\Ep(\bar{\eta})_f\big)
    \lra \bigoplus_{w\in\Si(F_\infty)}\frac{H^1\big(F_{n,w},\Ep(\bar{\eta})_f\big)}
    {C_w\big(\Ep(\bar{\eta})_f/F_{n}\big)}
        \Bigg),\]
        where \[ C_w\Big(\Ep(\bar{\eta})_f/F_{n}\Big) = \begin{cases}  \Big(E^{s_w}(\mF_{\infty,w})^{\eta}\ot\Qp/\Zp\Big)_f^{\Ga_n},& \mbox{if $w\in \Si_p^{ss}(F_{\infty})$}, \\
 \Bigg(\im\left(H^1(F_{n,w},\widehat{E}[p^{\infty}](\bar{\eta})_f)\lra H^1(F_{n,w}, E[p^{\infty}](\bar{\eta})_f)\right)\Bigg)_{div},& \mbox{if $w\in \Si_p^{ord}(F_\infty)$}, \\
 H^1_{ur}(F_{n,w},E[p^\infty](\bar{\eta})_f), & \mbox{if $w\nmid p$,} \end{cases}\]
and
\[ e_{\eta}\Sel^{\overrightarrow{s},str}(E[p^m]/\mF_n)_f\cong \ker\Bigg(H^1\big(G_{\Sigma}(F_n),E[p^m](\bar{\eta})_f\big)
    \lra \bigoplus_{w\in\Si(F_\infty)}\frac{H^1\big(F_{n,w},\Ep(\bar{\eta})_f\big)}
    {C_w\big(\Ep(\bar{\eta})_f/F_{n}\big)}
        \Bigg),\]
        where \[ C_w\Big(\Ep(\bar{\eta})_f/F_{n}\Big) = \begin{cases}  \Big(E^{s_w}(\mF_{n,w})^{\eta}\ot\Qp/\Zp\Big)^{\Ga_n}_f[p^m],& \mbox{if $w\in \Si_p^{ss}(F_{n})$}, \\
 \Bigg(\ker\left(H^1(F_{n,w},E[p^{m}](\bar{\eta})_f)\lra \displaystyle\frac{H^1(F_{n,w}, E[p^{\infty}](\bar{\eta})_f)}{C_w\big(\Ep(\bar{\eta})_f/F_{n}\big)}\right)\Bigg),& \mbox{if $w\in \Si_p^{ord}(F_n)$}, \\
 H^1_{ur}\big(F_{n,w}, E[p^{m}](\bar{\eta})_f\big), & \mbox{if $w\nmid p$.} \end{cases}\]
\el

\bpf
 By a similar argument to that in Lemma \ref{coprime H1}, we have the identification $e_{\eta}H^1(G_{\Sigma}(\mF_n),\Ep)_f\cong H^1\big(G_\Si(F_{n}),\Ep(\bar{\eta})_f\big)$. It therefore remains to show that the $e_{\eta}$-component of the local condition in the definition of the signed Selmer group $\Sel^{\s,str}(E/\mF_{\infty})$ coincides with $C_w\Big(\Ep(\bar{\eta})_f/F_{\infty}\Big)$. This follows from a direct verification noting Lemmas \ref{twist and eta}, \ref{description of eC} and \ref{orthogonal away from p lemma}.\epf

\br
We have set up our notation in a way so that for each $w\in \Si(F_{n})$, one has
\[ e_{\eta}\left(\bigoplus_{u\in \Si(\mF_n),u|w}C_u\big(\Ep/\mF_{n}\big)\right)= C_w\big(\Ep(\bar{\eta})_f/F_{n}\big) \]
and
\[ e_{\eta}\left(\bigoplus_{u\in \Si(\mF_n), u|w}C_u\big(E[p^m]/\mF_{n}\big)\right)= C_w\big(E[p^m](\bar{\eta})_f/F_{n}\big). \]\er

\bp \label{Poitou-Tate}
For a character $\eta$ of $G=\Gal(F(\mu_p)/F)$ and a power series $f\in \La$ which is not divisible by $p$, we have
\[ \frac{\Big|e_{\eta}\Sel^{\overrightarrow{s},str}(E/\mF_n)_f\Big|}
{\Big|e_{\bar{\eta}}\Sel^{\overrightarrow{s},str}(E/\mF_n)_{f^{\iota}}\Big|} =\prod_{w\in\Si^{ord}_p(F_n)}\frac{\big| E(F_{n,w})[p^\infty](\bar{\eta})_f/p^m\big| \big|\widehat{E}(F_{n,w})[p^{\infty}](\bar{\eta})_f/p^m\big|}{\big|E(F_{n,w})[p^m](\bar{\eta})_f\big|}. \]
Furthermore, the expression on the right is bounded independent of $m$ and $n$ (for a fixed $f$).
\ep

\bpf
Each term on the right is bounded independent of $m$ and $n$ (for a fixed $f$) by Lemma \ref{Imai finite}. But since $F_{\infty}/F$ is the cyclotomic $\Zp$-extension, there are only finite number of primes in $\Si^{ord}_p(F_{\infty})$. Hence the expression on the right is bounded independent of $m$ and $n$ (for a fixed $f$). It therefore remains to prove the equality as asserted in the proposition. In view of Corollary \ref{supersingular orthogonal}, Lemmas \ref{ordinary orthogonal} and \ref{orthogonal away from p lemma}, the local conditions $C_w\Big(E[p^m](\bar{\eta})_f/F_{n}\Big)$ and $C_w\Big(E[p^m](\eta)_{f^{\iota}}/F_{n}\Big)$ are orthogonal complement to each other. Thus, by a Poitou-Tate duality argument (for instance, see \cite[pp 125]{G89} or \cite[pp 89]{Kim08}), we have
\[\frac{\Big|e_{\eta}\Sel^{\overrightarrow{s},str}(E/\mF_n)_f\Big|}
{\Big|e_{\bar{\eta}}\Sel^{\overrightarrow{s},str}(E/\mF_n)_{f^{\iota}}\Big|} =
 \frac{\Big|H^1\big(G_{\Si}(F_n), E[p^m](\bar{\eta})_f\big)\Big|}{\Big|H^2\big(G_{\Si}(F_n), E[p^m](\bar{\eta})_f\big)\Big|} \prod_{w\in\Si^{ord}_p(F_n)} \frac{\Big|C(E[p^m](\bar{\eta})_f/F_n)\Big|\Big|H^2\big(F_n, E[p^m](\bar{\eta})_f\big)\Big|}{\Big|H^1\big(F_n, E[p^m](\bar{\eta})_f\big)\Big|},\]
 where in the calculation, we have made use of the fact that $H^0\big(G_{\Si}(F_n), E[p^m](\eta)_{f^{\iota}}\big) =0$ (which is a consequence of Lemma \ref{supersingular points twist}). Again, taking Lemma \ref{supersingular points twist} into account, the global Euler characteristic formula tells us that
 \[ \frac{\Big|H^1\big(G_{\Si}(F_n), E[p^m](\bar{\eta})_f\big)\Big|}{\Big|H^2\big(G_{\Si}(F_n), E[p^m](\bar{\eta})_f\big)\Big|} = p^{\big(r_1(F_n) +2r_2(F_n)\big)m\deg(f)} =  p^{[F_n:\Q]m\deg(f)},\]
 where $r_1(F_n)$ (resp., $r_2(F_n)$) is the number of real places of $F_n$ (resp., number of pairs of complex places of $F_n$). The required formula now follows from combining this calculation with Lemmas \ref{local ordinary calculation}, \ref{local supersingular lemma} and \ref{local away from p lemma}.
\epf

We can now prove the main theorem of the paper.

\bpf[Proof of Theorem \ref{main theorem}]
We shall verify that the hypotheses of Proposition \ref{compare La modules} are satisfied. By virtue of Proposition \ref{strict = usual}, we are reduced to showing that for a fixed $f$, the quantity  \[\frac{\left|\Big(e_{\eta}\Sel^{\overrightarrow{s},str}(E/\mF_\infty)_f\Big)^{\Ga_n}[p^m]\right|}
{\left|\Big(e_{\bar{\eta}}\Sel^{\overrightarrow{s},str}(E/\mF_\infty)_{f^{\iota}}\Big)^{\Ga_n}[p^m]\right|} \] is bounded independent of $m$ and $n$.
Since $e_{\eta}(-)_f$ is exact on sequence of $\La[G]$-modules, it follows from Propositions \ref{control theorem for strict} and \ref{control theorem for strict torsion} that
the kernel and cokernel of the maps
\[ e_{\eta}\Sel^{\s,str}(E/\mF_n)_f[p^m] \lra e_{\eta}\Sel^{\overrightarrow{s},str}(E/\mF_\infty)_f^{\Ga_n}[p^m]\]
\[ e_{\eta}\Sel^{\s,str}(E[p^m]/\mF_n)_f \lra e_{\eta}\Sel^{\overrightarrow{s},str}(E/\mF_n)_f[p^m]\]
are finite and bounded independent of $m$ and $n$. It follows from this that we are reduced to showing that
 the quantity  \[\frac{\left|e_{\eta}\Sel^{\s,str}(E[p^m]/\mF_n)_f\right|}
{\left| e_{\bar{\eta}}\Sel^{\s,str}(E[p^m]/\mF_n)_{f^{\iota}}\right|} \] is bounded independent of $m$ and $n$. But this is precisely the assertion of Proposition \ref{Poitou-Tate}. Hence we have completed the proof of our theorem.
\epf

\subsection{Proof of Corollary \ref{main corollary}} \label{proof of main corollary subsec}

We retain the setting of the previous subsections.
Recall that by Lemmas \ref{strict = usual} and \ref{twist selmer}, one has
\[ e_{\eta}\Sel^{\s}(E/\mF_{\infty}) \cong \ker\Bigg(H^1\big(G_{\Sigma}(F_\infty),\Ep(\bar{\eta})\big)
    \lra \bigoplus_{w\in\Si(F_\infty)}\frac{H^1\big(F_{\infty,w},\Ep(\bar{\eta})\big)}
    {C_w\big(\Ep(\bar{\eta})/F_{\infty}\big)}
        \Bigg),\]
        where \[ C_w\Big(\Ep(\bar{\eta})/F_{\infty}\Big) = \begin{cases}  E^{s_w}(\mF_{\infty,w})^{\eta}\ot\Qp/\Zp,& \mbox{if $w\in \Si_p^{ss}(F_{\infty})$}, \\
 \im\left(H^1(F_{\infty,w},\widehat{E}[p^{\infty}](\bar{\eta}))\lra H^1(F_{\infty,w}, E[p^{\infty}](\bar{\eta}))\right),& \mbox{if $w\in \Si_p^{ord}(F_\infty)$}, \\
 0, & \mbox{if $w\nmid p$.} \end{cases}\]

To continue, we require the following result which is a slight refinement of \cite[Theorem 3.14]{KimPM} and \cite[Theorem 1.3]{KO} (also see \cite[Proposition 4.14]{G99}).

\bp \label{nonexistence of finite}
If $e_{\eta}\Sel^{\s}(E/\mF_{\infty})$ is cotorsion over $\La$, then $\Big(e_{\eta}\Sel^{\s}(E/\mF_{\infty})\Big)^{\vee}$ has no nontrivial finite $\La$-submodules.
\ep

\bpf
 Observe that it follows from Lemma \ref{supersingular points twist} that one has $E(F_{\infty})[p^{\infty}](\bar{\eta})=0$. The conclusion of the proposition will now follow from a similar argument to that in \cite[Proposition 4.14]{G99} or \cite[Theorem 3.14]{KimPM}.
\epf

We can now give the proof of Corollary \ref{main corollary}.

\bpf[Proof of Corollary \ref{main corollary}]
The first assertion is immediate from the rank equality of Theorem \ref{main theorem}. We now proceed proving the second assertion. Suppose that $e_{\eta}X^{\s}(E/\mF_{\infty})=0$. Then by Theorem \ref{main theorem}, $e_{\bar{\eta}}X^{\s}(E/\mF_{\infty})^{\iota}$ is a finite $\La$-module, and in particular, a torsion $\La$-module. Thus, we may apply Proposition \ref{nonexistence of finite} to conclude that $e_{\bar{\eta}}X^{\s}(E/\mF_{\infty})$ has no nontrivial finite $\La$-submodule. Hence we must have $e_{\bar{\eta}}X^{\s}(E/\mF_{\infty})=0$.
\epf

\end{document}